\documentclass[12pt]{article}
\usepackage{amsfonts}

\usepackage{latexsym}
\usepackage{amssymb}
\usepackage{amsmath}
\usepackage{amsthm}
\usepackage[mathscr]{eucal}
\usepackage{graphicx}
\usepackage{hyperref}
\usepackage{caption}
\usepackage{subcaption}

\newtheorem{Lemma}{Lemma}

\newtheorem{Theorem}{Theorem}

\newtheorem{Definition}{Definition}

\renewcommand{\qed}{\hfill{\ \ \rule{2mm}{2mm}} \vspace{0.2in}}

\begin{document}

\title{Linearity Property of Unique Colourings in Random Graphs}
\author{ \textbf{Ghurumuruhan Ganesan}
\thanks{E-Mail: \texttt{gganesan82@gmail.com} } \\
\ \\
IISER Bhopal, India}
\date{}
\maketitle


\begin{abstract}
In this paper, we study unique colourings in random graphs as a generalization of both conflict-free and injective colourings. Specifically, we impose the condition that a \emph{fraction} of vertices in the neighbourhood of any vertex are assigned unique colours and use vertex partitioning and the probabilistic method to show that the minimum number of colours needed grows \emph{linearly} with the uniqueness parameter, unlike both conflict-free and injective colourings. We argue how the unboundedness of the vertex neighbourhoods influences the linearity property and illustrate our case with an example involving tree-unique colourings in random graphs.

\vspace{0.1in} \noindent \textbf{Key words:} Unique Colourings, Conflict-free Colourings, Injective Colourings, Linearity Property.

\vspace{0.1in} \noindent \textbf{AMS 2000 Subject Classification:} Primary: 05C62;
\end{abstract}

\bigskip

\renewcommand{\theequation}{\thesection.\arabic{equation}}
\setcounter{equation}{0}
\section{Introduction} \label{intro}
Conflict-free and injective colourings of random graphs are important problems from both theoretical and application perspectives. Conflict-free colouring (Pach and Tardos (2009)) of a graph requires that the closed neighbourhood of each vertex contains a unique colour and at the other end, injective colouring (Hahn et al. (2002)) implies that all the vertices in the open neighbourhood receive unique colours. For deterministic graphs, many bounds for various classes are known (see Pach and Tardos (2009), del R\'io-Chanona et al. (2016) and references therein). Recently, these colouring parameters have also been studied from the random graph context. Glebov et al. (2014) obtained sharp bounds for the conflict-free colouring number of random graphs and later del R\'io-Chanona et al. (2016) established upper and lower bounds of same order for injective colouring.

In this paper, we consider unique colourings as a generalization of both conflict-free and injective colourings by imposing the condition that a \emph{fraction} of the vertices in each neighbourhood receive unique colours. We use vertex partitioning and the probabilistic method to show that the minimum required number of colours grows linearly with the uniqueness parameter, unlike both conflict-free and injective colourings. We argue how unboundedness of the vertex neighbourhoods affects the linearity property and illustrate with an example of tree-unique colouring.

The paper is organized as follows. In Section~\ref{sec_conf_free}, we state and prove our main result regarding the linearity property of the colouring number for unique colourings. We also describe how the unbounded neighbourhood property influences the linear relationship of the colouring number with the uniqueness parameter and in Section~\ref{sec_bound}, illustrate our argument with an example of tree-unique colouring.

\renewcommand{\theequation}{\arabic{section}.\arabic{equation}}
\setcounter{equation}{0}
\section{Unique Colourings}\label{sec_conf_free}
For~\(n \geq 2\) let~\(K_n\) be the complete graph on~\(n\) vertices and let~\(\{Z(f)\}_{f \in K_n}\) be independent random variables indexed by the edge set of~\(K_n\) and with distribution
\begin{equation}\label{x_dist}
\mathbb{P}(Z(f) = 1) = p = 1-\mathbb{P}(Z(f) = 0)
\end{equation}
for some~\(0 < p  < 1.\) We let~\(G\) be the random graph formed by the union of all edges~\(f\) satisfying~\(Z(f) = 1.\) An edge with endvertices~\(u\) and~\(v\) is denoted as~\((u,v)\) and the degree of a vertex~\(u\) in~\(G\) is the number of edges in~\(G\) containing~\(u\) as an endvertex. Equivalently if~\(N(u)\) denotes the set of all vertices adjacent to~\(u\) in~\(G,\) then the degree~\(d(u)\) of~\(u\) equals~\(\#N(u),\) the cardinality of~\(N(u).\)

For integer~\(k \geq 1,\) a~\(k-\)colouring of~\(G\) is a map~\(f : \{1,2,\ldots,n\} \rightarrow \{1,2,\ldots,k\}.\) The colouring~\(f\) is said to be \emph{proper} if~\(f(u) \neq f(v)\) whenever the edge~\(e=(u,v)\) with endvertices~\(u\) and~\(v\) belongs to~\(G.\) Suppose~\(f\) is a~\(k-\)colouring (not necessarily proper) of~\(G.\) A vertex~\(v\) in the neighbourhood~\(N(u)\) of a vertex~\(u\) is said to be~\((f,N(u))-\)\emph{unique} or simply unique if~\(f(v) \neq f(w)\) for all~\(w \in N(u).\)
\begin{Definition}\label{conf_def}
For~\(0 < \eta < 1,\) we say that~\(f\) is a~\(\eta-\)\emph{injective} colouring of~\(G\) if for each vertex~\(u\) at least a fraction~\(\eta\) of the vertices in the neighbourhood~\(N(u)\) receive unique colours.

We denote~\(\chi_{i,\eta}(G)\) to be the minimum number of colours needed for an~\(\eta-\)injective colouring of~\(G.\)
\end{Definition}
For~\(\eta = 1\) the above definition reduces to the usual notion of injective colouring and it is well known (del R\'io-Chanona et al. (2016)) that~\(\chi_{i,1}(G) \geq \frac{C (np)^2}{\log{n}}\) with high probability, i.e. with probability converging to one as~\(n \rightarrow \infty\) for some constant~\(C  > 0,\) where, throughout, logarithms are natural and constants do not depend on~\(n.\)

For partial injective colourings as described in Definition~\ref{conf_def}, we have the following result regarding the corresponding colouring number.
\begin{Theorem}\label{cor_one}
For every~\(\gamma > 0\) and~\(0 < \eta < 1,\) there are positive constants~\(M = M(\gamma,\eta)\) and~\(D = D(\gamma,\eta)\) such that if~\(p \geq \frac{M\log{n}}{n}\)  then
\begin{equation}\label{thmix2}
\mathbb{P}\left(\frac{\eta np}{2} \leq \chi_{i,\eta}(G) \leq D np\right) \geq 1- \frac{1}{n^{\gamma}}.
\end{equation}
\end{Theorem}
In other words, for any fraction~\(0 < \eta < 1\) the corresponding partial injective colouring number grows \emph{linearly} in~\(np\) unlike the ``full" injective colouring with~\(\eta = 1,\) described in the previous paragraph.

We prove Theorem~\ref{cor_one} as a Corollary of the following Lemma of independent interest, regarding unique colourings. Formally, we say that a colouring~\(f\) of~\(G\) is~\(r-\)\emph{unique} if for each vertex~\(u,\) there are at least~\(\min(r,d(u)+1)\) unique colours in the \emph{closed} neighbourhood~\(N(u) \cup \{u\}.\)  We denote~\(\chi_{U,r}(G)\) to be the minimum number of colours needed for an~\(r-\)unique colouring of~\(G.\)

The case~\(r=1\) is also known as the conflict-free colouring (Glebov et al. (2014)) of~\(G\) and the case~\(r=\infty\) could be interpreted as ``strong" injective colouring. Clearly any proper colouring is conflict-free and any strong injective colouring is proper while the corresponding converses are not necessarily true. Consequently~\(\chi_{U,1}(G) \leq \chi(G) \leq \chi_{U,\infty}(G)\) where~\(\chi(G)\) denotes the minimum number of colours needed for a proper colouring of~\(G.\)  It is well-known (see for example, Ganesan (2020)) that~\(\chi(G)\) grows as~\(\frac{np}{\log{n}}\) with high probability and from (Glebov et al. (2014)) we get that~\(\chi_{U,1}(G)\) grows as~\(\log{n}\) with high probability. Of course, from the discussion following Definition~\ref{conf_def}, we know that~\(\chi_{U,\infty}(G)\) is at least of the order of~\(\frac{(np)^2}{\log{n}}\) with high probability.

The following result obtains bounds for the~\(r-\)unique colouring number of~\(G\) for general~\(r.\)
\begin{Lemma}\label{thm_a}  For every~\(\gamma > 0\) and~\(0 < \theta < 1\) there are positive constants~\(M_0 = M_0(\gamma,\theta)\) and~\(D_0 = D_0(\gamma,\theta)\) such that if~\(p \geq \frac{M_0\log{n}}{n}\)  and~\(1 \leq r \leq \theta np,\) then
\begin{equation}\label{thmix}
\mathbb{P}\left(r \leq \chi_{U,r}(G) \leq D_0 \max\left(r, \log{n}\right)\right) \geq 1- \frac{1}{n^{\gamma}}.
\end{equation}
\end{Lemma}
In particular for~\(r \geq \log{n},\) we see that with high probability~\(\chi_{U,r}\) grows linearly in the uniqueness parameter~\(r.\) In fact, in our proof of Theorem~\ref{thm_a} below, we crucially use the fact that each vertex has roughly~\(np \longrightarrow \infty\) neighbours with high probability and employ the probabilistic method within these unbounded neighbourhoods to obtain the desired colouring using~\(O(r)\) colours. In contrast, even~\(1-\)unique colouring within \emph{bounded} sub-neighbourhoods may require of the order of~\(np\) colours and we highlight this aspect with an example in the next Section.




To prove Lemma~\ref{thm_a}, we use the following standard deviation estimate which we state separately for convenience:  Let~\(\{X_j\}_{1 \leq j \leq r}\) be independent Bernoulli random variables with~\(\mathbb{P}(X_j = 1) = 1-\mathbb{P}(X_j = 0) > 0.\) If~\(T_r := \sum_{j=1}^{r} X_j,\theta_r := \mathbb{E}T_r\) and~\(0 < \gamma \leq \frac{1}{2},\) then
\begin{equation}\label{conc_est_f}
\mathbb{P}\left(\left|T_r - \theta_r\right| \geq \theta_r \gamma \right) \leq 2\exp\left(-\frac{\gamma^2}{4}\theta_r\right)
\end{equation}
for all \(r \geq 1.\) For a proof of~(\ref{conc_est_f}), we refer to Corollary A.1.14, pp. 312 of Alon and Spencer (2008).

\emph{Proof of Lemma~\ref{thm_a}}: The lower bound follows in~(\ref{thmix}) follows directly from the definition. Suppose~\(p \geq \frac{M\log{n}}{n}\) for some constant~\(M \geq 2\) to be determined later. We first prove the upper bound in~(\ref{thmix}) for the case~\(r \leq \frac{np}{2M}.\) Let~\(V_i, 1 \leq i \leq t = M \max\left(r,\frac{\log{n}}{2}\right)+1\) be~\(t\) disjoint vertex sets, each containing~\(\frac{1}{p}\) vertices, assuming for simplicity that~\(\frac{1}{p}\) is an integer. Such a partition is possible since each of~\(\frac{Mr}{p},\frac{M\log{n}}{2p}\) and~\(\frac{1}{p}\) is at most~\(\frac{n}{2}.\)

Assign colour~\(i\) to all the vertices of the set~\(V_i\) and for the rest of the vertices outside~\(\{V_j\}_{1 \leq j \leq t},\) we assign the colour~\(t+1.\) The probability that a vertex~\(v \notin V_1\) contains a unique neighbour in~\(V_1\) is~\(\alpha := (1-p)^{\frac{1}{p}-1} \geq \frac{1}{2e}\) for all~\(n\) large, since~\(p = o(1)\) i.e.,~\(p = p(n) \longrightarrow 0\) as~\(n \rightarrow \infty.\) Therefore if~\(X_{uni} = X_{uni}(v)\) is the number of sets in~\(\{V_j\}_{1 \leq j \leq t}\) containing a unique neighbour of~\(v,\) then~\(X_{uni}\) is Binomially distributed with parameters~\(t\) and~\(\alpha\) and so~\(\mathbb{E}e^{-X_{uni}} = \left(1-\alpha + \alpha e^{-1}\right)^{t} \leq  \exp\left(-(1-e^{-1})t\alpha\right).\)
From the Chernoff bound we therefore have for~\(\theta = \frac{1-e^{-1}}{2}\) and~\(M \geq \frac{8e^2}{e-1}(\gamma+1)\) that
\begin{equation}
\mathbb{P}\left(X_{uni} \leq \theta t \alpha\right) \leq \exp\left(-\frac{(1-e^{-1})}{2}t\alpha\right) \leq \exp\left(-\frac{M(e-1)}{8e^2} \log{n}\right) \leq \frac{1}{n^{\gamma+1}},\label{cherru}
\end{equation}
where the second inequality in~(\ref{cherru}) is true since~\(\alpha \geq \frac{1}{2e}\) and~\(t \geq \frac{M\log{n}}{2}.\) We fix such an~\(M\) henceforth.

For a vertex~\(v \in V_j,\) we define~\(X_{uni}(v)\) to be the number of sets in~\(\{V_l\}_{l \neq j}\) containing a unique neighbour of~\(v.\) An analogous analysis as above yields~(\ref{cherru}) and so from the union bound we get that~\(\mathbb{P}\left(\bigcup_{v} \left\{X_{uni}(v) \leq \theta t \alpha\right\}\right) \leq \frac{1}{n^{\gamma}}.\)
In other words, with probability at least~\(1-\frac{1}{n^{\gamma}},\) each vertex has at least~\(\theta t \alpha\) unique neighbours in the sets~\(\{V_j\}_{1 \leq j \leq t}.\) Using the estimates~\(t \geq M r, \alpha \geq \frac{1}{2e}, \theta = \frac{e-1}{2e}\) and~\(M \geq \frac{8e^2}{e-1} > 1,\) we get that~\(\theta t \alpha \geq r\) and this obtains~(\ref{thmix}) for the case~\(r \leq \frac{np}{2M}.\)

For values of~\(r > \frac{np}{2M},\) we use the probabilistic method. First the expected degree of any vertex is~\((n-1)p\) and so if~\(d(v)\) is the degree of vertex~\(v\) in~\(G,\) then arguing as in~(\ref{cherru}), we see that for any constant~\(0 < \epsilon < 1\)  there exists a constant~\(\mu > 0\) such that
\[\mathbb{P}\left(np(1-\epsilon) \leq d(v) \leq np(1+\epsilon)\right) \geq 1-e^{-\mu np} \geq 1-\frac{1}{n^{\gamma+1}}\]
provided~\(p \geq \frac{M\log{n}}{n}\) and we choose~\(M\) larger if necessary. Fixing such an~\(M\) and letting~\(E_{deg} := \bigcap_{v} \{np(1-\epsilon) \leq d(v) \leq np(1+\epsilon)\},\) we then get from the union bound that
\begin{equation}\label{deg_est}
\mathbb{P}(E_{deg}) \geq 1-\frac{1}{n^{\gamma}}.
\end{equation}

Suppose~\(E_{deg}\) occurs and we colour the vertices of~\(G\) randomly with colours chosen from the set~\(\{1,2,\ldots,q\},\) where~\(q = Dnp\) and~\(D \geq 2\) is a constant to be chosen later. Formally, we let~\(\{X_i\}_{1 \leq i \leq n}\) be independent random variables uniformly chosen from~\(\{1,2,\ldots,q\}\) that are also independent of~\(G\) and assign colour~\(X_i\) to the vertex~\(i.\) We also let~\(\mathbb{P}_X\) denote the distribution of~\((X_1,\ldots,X_n).\)

We now estimate the~\(\mathbb{P}_X-\)probability that the vertex~\(v\) has at most~\(r\) unique colours in its neighbourhood~\(N(v).\) Let~\(d= d(v) = \#N(v) \leq 2np\) be the degree of the vertex~\(v.\) If~\(N(v)\) has~\(k \leq r\) unique colours, then removing the vertices with these colours from~\(N(v),\) we are left with~\(d-k\) vertices coloured using at most~\(\frac{d-k}{2}\) colours. The number of ways of choosing the~\(k\) distinct colours is~\({q \choose k} \) and the number of ways of assigning these colours to vertices in~\(N(v)\) is at most~\(d^{k}.\)


Similarly, from the remaining~\(q-k\) colours, the number of ways of~\(l \leq \frac{d-k}{2}\) colours is~\({q \choose l} \leq {q \choose \frac{d-k}{2}},\) by the monotonicity of the Binomial coefficient and the fact~\(d= d(v) \leq 2np \leq q.\) The number of ways of assigning these~\(l \leq \frac{d-k}{2}\) colours to the remaining~\(d-k\) vertices (allowing for multiplicity) is at most~\(\left(\frac{d-k}{2}\right)^{d-k}.\) Since each colour assignment to the vertices of~\(N(v)\) occurs with probability~\(\frac{1}{q^{d}}\) and~\(l \leq d,\) we see that the~\(\mathbb{P}_X-\)probability that there at most~\(r\) unique colours in~\(N(v)\) is at most~\(\zeta := \frac{d}{q^{d}} \cdot \sum_{k \leq r} {q \choose k} d^{k} \cdot {q \choose \frac{d-k}{2}} \left(\frac{d-k}{2}\right)^{d-k}.\)

Using~\({q \choose l} \leq \left(\frac{qe}{l}\right)^{l},\) we get that~\({q \choose \frac{d-k}{2}} \left(\frac{d-k}{2}\right)^{d-k}\) is at most~\[\left(\frac{2eq}{d-k}\right)^{\frac{d-k}{2}} \left(\frac{d-k}{2}\right)^{d-k} = \left(\frac{qe(d-k)}{2}\right)^{\frac{d-k}{2}} \leq (qed)^{\frac{d-k}{2}}\]
and so~\(\zeta \leq d \sum_{k \leq r} \frac{1}{q^{k}} {q \choose k}d^{k} \cdot  \left(\frac{ed}{q}\right)^{\frac{d-k}{2}} \leq  d \left(\frac{ed}{q}\right)^{\frac{d-r}{2}} \sum_{k \leq r} \frac{1}{q^{k}} {q \choose k}d^{k},\)
since~\(k \leq r.\) Further using~\({q \choose k} \leq \frac{q^{k}}{k!}\) we get that
\begin{equation}
\zeta \leq d \left(\frac{ed}{q}\right)^{\frac{d-r}{2}}  \sum_{k} \frac{d^{k}}{k!} \leq de^{d} \left(\frac{ed}{q}\right)^{\frac{d-r}{2}}. \label{stich}
\end{equation}

Because~\(E_{deg}\) occurs, we know that~\((1-\epsilon) np \leq d \leq 2np\) and by definition~\(k \leq r \leq \theta np\) and so using~\(q = D np,\) we get that~\(\left(\frac{ed}{q}\right)^{\frac{d-k}{2}} \leq \left(\frac{2e}{D}\right)^{\frac{np}{2}(1-\epsilon-\theta)}.\) Choosing~\(\epsilon > 0\) small so that~\(1-\theta-\epsilon > 0,\) the estimate~(\ref{stich}) gives us that~\(\zeta \leq (2np) e^{2np} \left(\frac{2e}{D}\right)^{\frac{np}{2}(1-\epsilon-\theta)} \leq e^{-np}\) for all~\(n\) large, provided the constant~\(D = D(\epsilon,\theta)\) is chosen large enough. Fixing such a~\(D,\) we see that vertex~\(v\) has at most~\(r\) unique colours in its neighbourhood with~\(\mathbb{P}_X-\)probability at most~\(e^{-np}\) and so by the union bound,~\((X_1,\ldots,X_n)\) is an~\(r-\)unique colouring with~\(\mathbb{P}_X-\)probability at least~\(1-ne^{-np} \geq 1-\frac{1}{n} >0,\) provided~\(np \geq 2\log{n}.\)

Summarizing, if~\(E_{deg}\) occurs, then there exists an~\(r-\)unique colouring of~\(G\) and so from~(\ref{deg_est}), we again get~(\ref{thmix}). This completes the proof of Lemma~\ref{thm_a}.~\(\qed\)


\emph{Proof of Theorem~\ref{cor_one}}: Given~\(0 < \eta < 1\) choose~\(0 < \epsilon < \frac{1}{2}\) small and~\(0 < \theta <1\) such that~\(\eta(1+\epsilon) < \theta < \eta(1+2\epsilon) <1.\) Denoting~\(E_{\theta}\) to be the event that there exists a~\(\theta np-\)unique colouring of~\(G\) using at most~\(D \theta np\) colours, the estimate~(\ref{thmix}) in Lemma~\ref{thm_a} implies that~\(\mathbb{P}(E_{\theta}) \geq 1-\frac{1}{n^{\gamma}},\) provided~\(p \geq \frac{M\log{n}}{n}\) and~\(M > 0\) is large. Fixing such an~\(M,\) we now let~\(E_{deg} = E_{deg}(\epsilon)\) be the event as defined prior to~(\ref{deg_est}) so that each vertex has a degree between~\((1-\epsilon)np\) and~\((1+\epsilon)np.\) Choosing~\(M\) larger if necessarily, we again get from~(\ref{deg_est}) that~\(\mathbb{P}(E_{deg}) \geq 1-\frac{1}{n^{\gamma}}.\) Consequently, the union bound implies that\\\(\mathbb{P}(E_{deg} \cap E_{\theta}) \geq 1-\frac{2}{n^{\gamma}}.\)

If~\(E_{deg} \cap E_{\theta}\) occurs, then~\(\chi_{\eta}(G) \geq \eta\min_{u} d(u) \geq \eta (1-\epsilon) np \geq \frac{\eta np}{2}\) since\\\(\epsilon < \frac{1}{2}\) and moreover using~\(\max_{u}d(u) \leq (1+\epsilon)np,\) we get from our choice of~\(\theta\) that~\(\chi_{\eta}(G) \leq D_0 \theta np \leq D_0 \eta(1+2\epsilon) np \leq 2D_0\eta np\) for some constant\\\(D_0 = D_0(\gamma,\theta) > 0,\) again since~\(\epsilon < \frac{1}{2}.\)~\(\qed\)

\renewcommand{\theequation}{\arabic{section}.\arabic{equation}}
\setcounter{equation}{0}
\section{Bounded Sub-Neighbourhoods}\label{sec_bound}
Let~\(G = G(n,p)\) be the random graph with edge probability~\(p\) as defined in~(\ref{x_dist}). In this section, we illustrate the effect of neighbourhood size by considering unique colourings within \emph{bounded} sub-neighbourhoods of~\(G.\)



We begin with some definitions. For integer~\(k \geq 1\) let~\(f\) be a~\(k-\)colouring of~\(G,\) not necessarily proper. For a constant integer~\(t \geq 2\) let~\(T\) be any deterministic labelled tree on~\(t\) vertices and with vertex set~\(V(T).\) We say that a subgraph~\(H\) of~\(G\) with vertex set~\(V(H)\) is \emph{isomorphic} to~\(T\) if there exists a bijection~\(g : V(H) \rightarrow V(T)\) such that any two vertices~\(u\) and~\(v\) in~\(H\) are adjacent if and only if~\(g(u)\) and~\(g(v)\) are adjacent in~\(T.\)
\begin{Definition}\label{def_two} We say that~\(f\) is a~\(T-\)unique colouring of~\(G\) if every subgraph~\(H\) of~\(G\) isomorphic to~\(T\) contains at least one unique colour.

We let~\(\chi_T(G)\) denote the minimum number of colours needed for a~\(T-\)unique colouring of~\(G.\)
\end{Definition}
If~\(T\) is simply the single edge on~\(t=2\) vertices, then Definition~\ref{def_two} reduces to the usual notion of proper colouring described prior to Definition~\ref{conf_def}. 

For general~\(T,\) we have the following result regarding bounds for~\(\chi_T(G).\)
\begin{Theorem}\label{thm_b} If~\(p = o(1)\)  and~\(np \longrightarrow \infty\) as~\(n \rightarrow \infty,\) then
\begin{equation}\label{thulp}
\mathbb{P}\left(\frac{D_1np}{\log{n}} \leq \chi_T(G) \leq D_2 (np)^{2\left(1-\frac{1}{t}\right)}\right) \geq 1- e^{-D_3(\log{n})^2} -ne^{-D_3np}
\end{equation}
for some positive constants~\(D_i, 1 \leq i \leq 3\) and all~\(n\) large.
\end{Theorem}
In particular, allowing~\(T\) to be the star graph on~\(t\) vertices, we see from Theorem~\ref{thm_b} that even~\(1-\)unique colouring with \emph{bounded} sub-neighbourhood constraints requires at least of the order of~\(\frac{np}{\log{n}}\) colours unlike the ``unbounded" case which requires at most of the order of~\(\log{n}\) colours with high probability, as seen by plugging~\(r=1\) in Lemma~\ref{thm_a}.

To prove the upper bound in Theorem~\ref{thm_b} we use the following version of the Lov\'asz Local Lemma (Lemma~\(5.1.1.\) pp. 64, Alon and Spencer (2008)).
\begin{Lemma}\label{local} Let~\(A_1,\ldots,A_t\) be events in an arbitrary probability space. Let~\(\Gamma \) be the dependency graph for the events~\(\{A_i\},\) with vertex set~\(\{1,2,\ldots,t\}\) and edge set~\({\cal E};\) i.e. assume that each~\(A_i\) is independent of the family of events~\(A_j, (i,j) \notin {\cal E}.\) If there are reals~\(0 \leq y(i) < 1\) such that\\\(\mathbb{P}(A_i) \leq y(i) \prod_{(i,j)\in {\cal E}} (1-y(j))\) for each~\(i,\) then~\[\mathbb{P}\left(\bigcap_{i} A^c_i\right) \geq \prod_{1 \leq i \leq n} (1-y(i)) > 0.\]
\end{Lemma}

To prove the lower bound in Theorem~\ref{thm_b}, we use the following Lemma of independent interest regarding the probability of the event~\(E_T\) that~\(G\) contains a subgraph isomorphic to~\(T.\)
\begin{Lemma}\label{t_lem} If~\(p =o(1)\) and~\(np \longrightarrow \infty\) as~\(n \rightarrow \infty,\) then
\begin{equation}\label{thmulp}
\mathbb{P}(E_T) \geq 1-e^{-Cn^2p}
\end{equation}
for some constant~\(C = C(t) > 0\) depending only on~\(t\) and not on the choice of~\(T.\)
\end{Lemma}
We recall from Corollary~\(4.6\) pp.~\(82\) of Bollob\'as (2001) that if~\(p \geq n^{-\frac{t}{t-1}},\) then~\(E_T\) occurs with high probability. We use the estimate in~(\ref{thmulp}) obtained under the stronger assumption that~\(np \longrightarrow \infty,\) together with the union bound to obtain the lower bound for~\(\chi_T\) in the proof of Theorem~\ref{thm_b} below.

\begin{figure}[tbp]
\centering
\includegraphics[width=3.5in, trim= 20 400 20 270, clip=true]{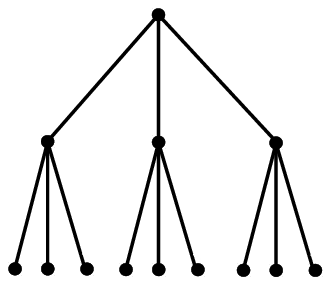}
\caption{The rooted tree~\(T_{3,2}.\)}
\label{fig_tree}
\end{figure}

\emph{Proof of Lemma~\ref{t_lem}}: We actually estimate the probability that the rooted ``nearly" regular tree~\(T_{l,t}\) with root degree~\(t\) and common vertex degree~\(t+1\)  and having~\(l\) levels occurs as a subgraph of~\(G.\) In Figure~\ref{fig_tree}, we have illustrated~\(T_{3,2}\) as an example.

The first step in our proof is to find a vertex whose degree is of the order of~\(np\) with ``very" high probability. Indeed if~\(m\) denotes the number of edges of~\(G,\) then~\(\mathbb{E}m={n \choose 2}p \leq \frac{n^2p}{2}\) and so if~\(E_{up} := \left\{ m \leq n^2p\right\},\) then by the deviation estimate~(\ref{conc_est_f}) we have that
\begin{equation}\label{edge_est}
\mathbb{P}\left(E_{up}\right) \geq 1-e^{-C_1n^2p}
\end{equation}
for some constant~\(C_1 > 0.\) Next for a deterministic set~\(S\) containing~\(s \geq \frac{2n}{3}\) vertices, let~\(m(S)\) be the number of edges of~\(G\) containing both endvertices in~\(S.\) We have that~\(\mathbb{E}m(S) \geq {s \choose 2}p \geq \frac{s^2p}{4} \geq \frac{n^2p}{9}\) and so again from~(\ref{conc_est_f}), we get that~\(\mathbb{P}\left(m(S) \geq \frac{n^2p}{10} \right) \geq 1-e^{-2C_2n^2p}\) for some constant~\(C_2 > 0.\) Let~\(E_{low} := \bigcap_{S} \left\{m(S) \geq \frac{n^2p}{10}\right\}\) where the intersection is over all sets~\(S\) of size~\(s \geq \frac{2n}{3}.\) The number of choices for~\(S\) is at most~\(2^{n}\) and so by the union bound, we get that
\begin{equation}\label{edge_est2}
\mathbb{P}(E_{low}) \geq 1-2^{n}e^{-2C_2n^2p} \geq 1-e^{-C_2n^2p}
\end{equation}
for all~\(n\) large, since~\(np \longrightarrow \infty\) by Lemma statment.

Suppose~\(E_{low} \cap E_{up}\) occurs so that the total number of edges in~\(G\) is at most~\(n^2p.\) This in turn implies that the sum of vertex degrees is no more than~\(2n^2p\) and so at most~\(\frac{n}{3}\) vertices have degree larger than~\(6np.\) Equivalently, at least~\(\frac{2n}{3}\) vertices of~\(G\) have degree less than~\(6np.\) Since~\(E_{low}\) also occurs, at least one among these~\(\frac{2n}{3}\) vertices have degree at least~\(\frac{np}{10}.\) In effect, we have shown that if~\(E_{low} \cap E_{up}\) occurs, then there exists a vertex~\(z\) whose degree satisfies~\(\frac{np}{10} \leq d(z) \leq 6np.\)


If~\(E_{l}\) denotes the event that no subgraph of~\(G\) is isomorphic to~\(T_{l,t},\) then from the above discussion and~(\ref{edge_est}) and~(\ref{edge_est2}), we see that
\begin{eqnarray}
\mathbb{P}(E_l^c) &\leq& \mathbb{P}(E_l^c \cap E_{low} \cap E_{up}) + \mathbb{P}(E^c_{low} \cup E^c_{up}) \nonumber\\
&\leq&  \mathbb{P}(E_l^c \cap E_{low} \cap E_{up}) + e^{-C_1n^2p}  +e^{-C_2n^2p}\nonumber\\
&\leq& \mathbb{P}\left(E_l^c \bigcap \bigcup_{1 \leq z \leq n}\left\{\frac{np}{10} \leq d(z) \leq 6np\right\}\right) + e^{-C_1n^2p} +e^{-C_2n^2p}\nonumber\\
&\leq& \sum_{z=1}^{n} \mathbb{P}\left(E_l^c \bigcap \left\{\frac{np}{10} \leq d(z) \leq 6np\right\}\right) + e^{-C_1n^2p} +e^{-C_2n^2p}\nonumber\\
&=& n\mathbb{P}\left(E_l^c \bigcap \left\{\frac{np}{10} \leq d(1) \leq 6np\right\}\right) + e^{-C_1n^2p} +e^{-C_2n^2p} \label{thee}
\end{eqnarray}

Suppose~\(E_l^c \bigcap \left\{\frac{np}{10} \leq d(1) \leq 6np\right\}\) occurs and for a deterministic set~\(S\) of size between~\(\frac{np}{10}\) and~\(6np,\) let~\(N(1)=S\) be the set of neighbours of the vertex~\(1.\) The number of vertices in~\(S\) is at least~\(\frac{np}{10} \geq t\) since~\(np \rightarrow \infty\) (see Lemma statement) and similarly, the number of vertices outside~\(S\) is at least~\(n-6np \geq \frac{n}{2}\) since~\(p = o(1)\) again by Lemma statement. We now pick~\(t\) neighbours~\(s_1,\ldots,s_t\) from~\(S\) and divide~\(\frac{n}{2}\) of the vertices outside~\(S\) (picked according to a deterministic rule) into~\(t\) disjoint sets~\({\cal T}_j, 1 \leq j \leq t\) each of size~\(n_1 := \frac{n}{2t},\) assumed to be an integer for simplicity.

Let~\(G_j\) be the subgraph of~\(G\) induced by the vertices of~\({\cal T}_j \cup \{s_j\}.\) If~\(E_l^c\) occurs, then there must exist a subgraph~\(G_j ,1 \leq j \leq t\) that does not have~\(T_{l-1,t}\) as a subgraph. Therefore if~\(E_{l,j}\) denotes the event that~\(G_j\) contains a subgraph isomorphic to~\(T_{l-1,t}\) then
\begin{eqnarray}
\mathbb{P}\left(E_l^c \bigcap \left\{\frac{np}{10} \leq d(1) \leq 6np\right\}\right)
&=& \sum_{S} \mathbb{P}\left(E_l^c \bigcap \{N(1) = S\}\right) \nonumber\\
&\leq& \sum_{S} \mathbb{P}\left(\{N(1) = S\} \bigcap\bigcup_{1 \leq j \leq t} E_{l,j}^c\right) \nonumber\\
&\leq& t\sum_{S}\mathbb{P}\left(\{N(1) = S\} \bigcap E_{1,j}^c\right) \nonumber\\
&=& t\sum_{S}\mathbb{P}(N(1) = S) \mathbb{P}(E_{1,j}^c) \nonumber\\
&\leq& t \mathbb{P}(E_{1,j}^c) \label{thee2}
\end{eqnarray}
where the second equality in~(\ref{thee2}) follows from the fact that the event\\\(N(1) = S\) depends on the set of edges containing~\(1\) as an endvertex and is therefore independent of the event~\(E_{1,j}^c\) that depends only on edges \emph{not} containing~\(1\) as an endvertex.

Denoting~\(q_l(n) := \mathbb{P}(E_l^c)\) and combining~(\ref{thee}) and~(\ref{thee2}) and setting~\(D_1 := \max(C_1,C_2),\) we get the recursion~\(q_l(n) \leq tq_{l-1}\left(\frac{n}{2t}\right) + 2e^{-D_1n^2p}.\) Proceeding recursively, we get that
\begin{eqnarray}
q_l(n) &\leq& t \left(t q_{l-2}\left(\frac{n}{(2t)^2}\right) + 2e^{-D_2n^2p}\right) + 2e^{-D_1n^2p} \nonumber\\
&=& t^2q_{l-2}\left(\frac{n}{(2t)^2}\right) + 2te^{-D_2n^2p} + 2e^{-D_1n^2p} \nonumber\\
&\leq& \ldots \nonumber\\
&=& t^{l-1}q_1\left(\frac{n}{(2t)^{l-1}}\right) + \sum_{j=1}^{l-1} 2t^{j-1} e^{-D_jn^2p} \label{rec_eq2}
\end{eqnarray}
for some constants~\(D_j, 1 \leq j \leq l-1.\) The term~\(q_1(z)\) is simply  the probability that there is no vertex of degree at least~\(t\) in the random graph~\(G(z,p).\) Therefore arguing as in the paragraph containing~(\ref{edge_est}) and using~(\ref{edge_est}), we get that~\(q_1(z) \leq e^{-D_lz^2p}\) for some constant~\(D_l > 0\) and so from~(\ref{rec_eq2}), we get that~\(q_l(n) \leq 2\sum_{j=1}^{l}t^{j-1}e^{-D_jn^2p} \leq e^{-Dn^2p}\) for some constant~\(D > 0\) since~\(np \longrightarrow \infty.\)~\(\qed\)

\emph{Proof of Theorem~\ref{thm_b}}: We begin with the proof of the upper bound. The expected degree of any vertex in~\(G\) is~\((n-1)p\) and so if~\(E_{deg}\) denotes the event that each vertex in~\(G\) has degree at most~\(2np,\) then the deviation estimate~(\ref{conc_est_f}) together with the union bound implies that
\begin{equation}\label{deg_estax}
\mathbb{P}(E_{deg}) \geq 1-ne^{-C_1np}
\end{equation}
for some constant~\(C_1 >0.\)

Suppose~\(E_{deg}\) occurs and let~\(X_1,\ldots,X_n\) be independent random variables (also independent of~\(G\)) that are uniformly chosen from the set~\(\{1,2,\ldots,q\}\) for some integer~\(q \geq 1\) to be determined later. Assign colour~\(X_i\) to vertex~\(i\) and let~\(\mathbb{P}_X\) denote the distribution of~\((X_1,\ldots,X_n).\) Let~\({\cal T}\) be the set of all subgraphs of~\(G\) isomorphic to~\(T\) and for a tree~\(\tau \in {\cal T},\) let~\(A_{\tau}\) be the event that~\(\tau\) has a unique colour.

Recalling that~\(t\) denotes the number of vertices in~\(T\) and assuming that~\(t= 2z\)  is even for simplicity, we now estimate the~\(\mathbb{P}_X-\)probability that~\(\tau\) does not have a unique colour. Indeed, if the tree~\(\tau\) containing~\(2z\) vertices does not have a unique colour, then at most~\(z\) colours have been used to colour the vertices of~\(\tau.\) The number of ways of choosing~\(l \leq z\) colours is~\({q \choose l}\) and the number of ways of assigning these~\(l\) colours to the vertices of~\(\tau\) is at most~\(l^{2z}.\) Therefore if~\(q \geq 2z\) then~\(\mathbb{P}_X(A^c_{\tau}) \leq \sum_{1 \leq l \leq z} \frac{{q \choose l} l^{2z}}{q^{2z}} \leq \frac{z {q \choose z} z^{2z}}{q^{2z}}\)
due to the monotonicity of the binomial coefficient and the fact that~\(q \geq 2z.\) Further using the estimate~\({q \choose z}  \leq \left(\frac{q e}{z}\right)^{z},\) we then get that
\begin{equation}\label{px_est}
\mathbb{P}_X(A^c_{\tau}) \leq z \left(\frac{q e}{z}\right)^{z} \cdot \frac{z^{2z}}{q^{2z}} = z \left(\frac{ze}{q}\right)^{z} := \frac{y(\tau)}{2}.
\end{equation}

By definition, the events~\(A_{\tau}\) and~\(A_{\lambda}\) are dependent only if the trees~\(\tau\) and~\(\lambda\) share a common vertex, which we denote as~\(\tau \cap \lambda \neq \emptyset.\) Since~\(E_{deg}\) occurs, the degree of each vertex is at most~\(2np\) and so the number of subgraphs of~\(G\) isomorphic to~\(T\) and containing the vertex~\(1\) is at most~\(C_2(np)^{t-1} = C_2(np)^{2z-1}\) for some large constant~\(C_2 > 0.\) This in turn implies that~\(\tau \cap \lambda \neq \emptyset\) for at most~\(C_2(np)^{2z-1}\) trees~\(\lambda \in {\cal T}\) and so from~(\ref{px_est}), we see that there is a constant~\(C_3 > 0\) such that
\begin{eqnarray}
y(\tau) \prod_{\lambda : \lambda \cap \tau \neq \emptyset} (1-y(\lambda)) &\geq& 2z\left(\frac{ze}{q}\right)^{z} \left(1-2z\left(\frac{ze}{q}\right)^{z} \right)^{C_2(np)^{2z-1}} \nonumber\\
&\geq& 2z \left(\frac{ze}{q}\right)^{z} \left(1-C_3 \frac{(np)^{2z-1}}{q^{z}}\right)  \nonumber\\
&\geq& z\left(\frac{ze}{q}\right)^{z} \nonumber\\
&\geq& \mathbb{P}_X(A_{\tau}), \nonumber
\end{eqnarray}
provided~\(q^{z} \geq 2C_3 (np)^{2z-1}\) or equivalently~\(q \geq C_4 (np)^{2-\frac{1}{z}}\) for some constant~\(C_4 > 0.\) Fixing such a~\(q,\) the local lemma (Lemma~\ref{local}) then implies that there is a~\(T-\)unique colouring of~\(G\) with positive~\(\mathbb{P}_X-\)probability. From~(\ref{deg_estax}), we then obtain the upper bound in~(\ref{thulp}).

For the lower bound, we use Lemma~\ref{t_lem}. Let~\(C\) be the constant in Lemma~\ref{t_lem} and for a deterministic set~\(S\) containing~\(s := \frac{2\log{n}}{Cp}\) vertices, let~\(E_{tree}(S)\) be the event that the subgraph~\(G_S\) of~\(G\) induced by the vertices of~\(S,\) itself contains a subgraph isomorphic to~\(T.\) From Lemma~\ref{t_lem}, we see that~\(\mathbb{P}(E^c_{tree}(S)) \leq e^{-Cs^2p}\) and so if~\(E_{tree} := \bigcap_{S} E_{tree}(S)\) where the intersection is over all sets~\(S\) containing~\(s\) vertices, the union bound implies that~\(\mathbb{P}(E^c_{tree}) \leq {n \choose s}e^{-Cs^2p}.\) Using~\({n \choose s} \leq n^{s} = e^{s\log{n}}\) we then get that
\begin{equation}\label{milf}
\mathbb{P}(E^c_{tree}) \leq e^{s\log{n} - Cs^2p} = e^{-s\log{n}} \leq e^{-C_5(\log{n})^2}
\end{equation}
for some constant~\(C_5 > 0,\) since~\(s = \frac{2\log{n}}{Cp} \geq \frac{2\log{n}}{C}.\) If~\(E_{tree}\) occurs and~\(f\) is a~\(T-\)unique colouring of~\(G,\) then each colour is used on at most~\(s\) vertices and therefore~\(\chi_{T}(G) \geq \frac{n}{s} \geq \frac{Cnp}{2\log{n}}.\) The estimate~(\ref{milf}) then obtains the lower bound in~(\ref{thulp}).~\(\qed\)

\underline{\emph{Acknowledgement}}: I thank IMSc faculty for crucial comments and also thank IMSc and IISER Bhopal for my fellowships.


\bibliographystyle{plain}

\end{document}